\documentclass[12pt]{article}
\usepackage{amsfonts,amssymb,amscd} 
\usepackage{epsfig}
\usepackage{graphicx}
\usepackage{tikz-cd}
\usepackage[all]{xy}
\usepackage[T1]{fontenc}
\usepackage[utf8]{inputenc}
\begin{document}
\title{Free groups, covering spaces and Artin's theorem}
\author{
Gopala Krishna Srinivasan
\footnote{e-mail address: gopal@math.iitb.ac.in}
}
\date{
 Department of  Mathematics, Indian Institute of Technology Bombay, 
Mumbai 400 076
}
\maketitle
\section{Introduction:} 
One of the charming experiences from a first encounter with algebraic topology, in particular the theory of 
covering spaces, is the topological proofs of several important theorems in group theory. 
Massey [4] devotes his last chapter to some of these as a 
d\'enouement of the theory of fundamental groups and covering spaces.  
The theorem alluded to in the title of this note is the following  unpublished result of Emil Artin. 
\paragraph*{Theorem 1 (Artin:)} The commutator subgroup of the free group on two generators is not finitely generated.

This is a result that surely springs a  surprise and  belongs to  the  
family of results  in combinatorial group theory for which topological proofs as well as 
algebraic proofs are available. The latter usually employs the method of 
 Schreier transversals (see \cite{robinson}, p. 163). The book of John Stillwell \cite{stillwell}
is an excellent  account describing the close interactions between combinatorial group theory and topology. 
An authoritative historical account is the book of B. Chandler and W. Magnus \cite{chandler}
(see particularly pages 96-97).  John Stillwell \cite{stillwell} 
attributes the result to Emil Artin and provides a topological 
proof on page 101 
employing the infinite grid   
$$
G = \{ (x, y) \in \mathbb R^2 \;: \; x \in \mathbb Z \mbox{ or } y \in \mathbb Z\}
$$
as a covering space (the universal abelian cover) for the wedge of two circles $S^1\vee S^1$. 
A proof is also available in \cite{cohen} (p. 175).   
The present note provides yet another very short and crisp topological argument, also 
employing the infinite grid but takes a somewhat different route than the above mentioned proofs which the 
author believes would be of some value. 
As a further simplification we avoid in this note the  use of infinite graphs and 
their spanning trees. 
\section{Free groups and covering spaces:} 
For the benefit of the reader let us recapitulate some of the notions used in this note. If $V$ is a vector space, 
a subset $S \subset V$ is a basis of $V$ if the following two properties are satisfied: 
\begin{enumerate}
\item[(i)] For every vector space $W$, every map $f : S \longrightarrow W$ extends as  a linear transformation 
$F : V \longrightarrow W$ namely, the following diagram commutes: 
\[
\xymatrix{
S \ar[rr]^{\iota}\ar[rd]_{f} & & V \ar[ld]^{F}\\
& W
}
\]
where, $\iota : S \longrightarrow V$ is the inclusion map. 
\item[(ii)] 
The linear extension $F$ in (i) is unique. 
\end{enumerate} 
It is easy to see that this definition is equivalent to the standard definition given in elementary linear algebra courses. 
Indeed if $S$ fails to be linearly independent, the linear extension  $F$ in (i) may not exist while if $S$ fails to generate 
$V$ then the  map $F$ in (i) (if it exists) would not be unique. The above definition though more abstract, 
has the advantage that it readily generalizes 
yielding the notion of a  free group and a basis for a free group. 
\paragraph*{Definition 1:} A group $G$ is said to be free with basis $S \subset G$ if 
\begin{enumerate}
\item[(i)] For every group $H$, every map $f : S \longrightarrow H$ extends as  a  group homomorphism  
$F : G \longrightarrow H$ namely,  $F \circ i = f$. Again with $\iota : S \longrightarrow G$ denoting the inclusion, 
the following diagram commutes: 
\[
\xymatrix{
S \ar[rr]^{\iota}\ar[rd]_{f} & & G \ar[ld]^{F}\\
& H
}
\]
\item[(ii)] 
The extension $F$ in (i) is unique. 
\end{enumerate} 
\paragraph*{Remark:} If we replace every occurrence of the word ``group'' by ``abelian group'' we get the notion of a 
free abelian group $G$ with basis $S$. In the context of free abelian groups we have:
\paragraph*{Theorem 2:} Suppose $G$ is a free abelian group then, 
\begin{enumerate}
\item[(i)]
Any two bases of $G$ have the same cardinality and the common value is known as 
the rank of the free abelian group. 
\item[(ii)]
A subgroup of a free abelian group is again free abelian. 
\item[(iii)] If $G$ and $H$ are free abelian and $G \longrightarrow H$ is a surjective group homomorphism then rank of 
$H$ cannot exceed the rank of $G$. 
\end{enumerate}
\paragraph*{Proof:} See \cite{robinson} (pp. 100-101) for details on (i) and (ii). We remark that the proof proceeds by 
considering $G \otimes \mathbb Q$ as a  $\mathbb Q$ vector space. The map $g \mapsto g \otimes 1$ is a 
monomorphism of $G$ into $G \otimes \mathbb Q$ and 
if $S$ is a basis of $G$ then the set 
$S^* = \{s\otimes 1\;:\; s \in S\}$ 
is a basis of the vector space $G \otimes \mathbb Q$. From this one easily deduces (iii) or  use 
exercise 7(i) on p. 105 of \cite{robinson}. 
\hfill$\square$

If $G$ is any group denote by $[G, G]$ its commutator subgroup and by $A(G)$ its abelianization given by 
$A(G) = G/[G, G]$ and $\eta : G \longrightarrow A(G)$ the quotient map. We now prove that $A(G)$ is a 
free abelian group if 
$G$ is a free group. 
\paragraph*{Theorem 3:} Suppose $G$ is a free group with basis $S$, then 
\begin{enumerate}
 \item[(i)]
 $A(G)$ is a free abelian group with basis $\eta(S)$. 
 \item[(ii)] The restriction $\eta\vert_{S}$ is injective.
 \item[(iii)] Any two bases of $G$ have the same cardinality. The common value is called the rank of the free group. 
\end{enumerate}
\paragraph*{Proof:} From (i) and (ii) we immediately get (iii) upon invoking theorem 2. One can also 
give a direct proof \cite{robinson} (p. 50, exercise 7). 
Let us quickly dispose off (ii). 
If $\eta(s_1) = \eta(s_2)$ ($s_1 \neq s_2$) 
then $s_1^{-1}s_2$ lies in the commutator subgroup $[G, G]$. We take $H = \mathbb Z$ and declare 
$f(s_1) = 1$, $f(s) = 0$ if $s \neq s_1$ and $s \in S$. Then $f : S \longrightarrow \mathbb Z$  
extends as a group homomorphism $F : G \longrightarrow \mathbb Z$ 
and the target group being abelian (regarded additively), 
$F$ must map $[G, G]$ to the trivial element $\{0\}$.  Hence $F(s_1^{-1}s_2) = 0$ or $F(s_1) = F(s_2)$, which is a 
contradiction since $F(s_1) = f(s_1) = 1$ and $F(s_2) = f(s_2) = 0$.

Turning now to the proof of (i), let $H$ be an abelian group and $f: \eta(S) \longrightarrow H$ be any map. To show that 
$f$ extends to a group homomorphism $F : A(G) \longrightarrow H$, consider the map $f \circ \eta : S \longrightarrow H$ 
which must extend as a group homomorphism $\phi : G \longrightarrow H$. 
%$$
%Diagram
%$$
Now since $H$ is abelian, 
by the universal property of quotient there is a unique group homomorphism 
$F : A(G) \longrightarrow H$ such that $F \circ \eta = \phi$. Evaluating at a typical element of 
$S$ we get 
$$
F(\eta(s)) = \phi(s) = f(\eta(s)) 
$$
showing that $F$ indeed does the job. \hfill$\square$

From the theory of covering spaces we need theorem 4 below. 
It is well known that if $q : (\tilde{X}, \tilde x_0) \longrightarrow (X, x_0)$ is a covering projection the induced 
map $q_{*} : \pi_1(\tilde X, \tilde x_0) \longrightarrow \pi_1(X, x_0)$ is a monomorphism and the 
index of the subgroup $q_*(\pi_1(\tilde X, \tilde x_0))$ equals the cardinality of the fiber $q^{-1}(x_0)$. 
\paragraph*{Theorem 4:} Let $q : (\tilde{X}, \tilde x_0) \longrightarrow (X, x_0)$ be a covering projection between 
connected, locally path connected spaces. Then the following are equivalent:

(i) $q_*\pi_1(\tilde X, \tilde x_0)$ is a normal subgroup of $\pi_1(X, x_0)$.  

(ii) The group $\mbox{Deck}(\tilde{X}, X)$ of deck transformations acts transitively on the fibers. 

(iii) If $\gamma$ is a loop in $X$ based at $x_0$ then either every lift of $\gamma$ is closed or none of the 

lifts of $\gamma$ is closed. 
\newline
Further, when these conditions hold we have the following important relation:  
$$
\mbox{Deck}(\tilde{X}, X) = \pi_1(X, x_0)/q_*(\pi_1(\tilde X, \tilde x_0)) \eqno(1)
$$
\paragraph*{Proof:} See \cite{lima} (p. 153 and p. 167). \hfill$\square$ 
\paragraph*{Definition 2:}  A covering projection satisfying any of the above equivalent conditions 
is called a regular covering.

We shall now consider  the covering map 
$p : G \longrightarrow S^1\vee S^1$ is given by 
$$
p(x, y) = (\exp(2\pi i x), \exp(2\pi iy)). \eqno(2)
$$
This covering appears in [6] and [8] to illustrate that the fundamental group of 
$S^1\vee S^1$ is non-abelian. It is clear that the deck transformations of (2) are the translations:
$$
T_{(m, n)}(x, y) = (x+m, y+n),\quad m, n \in \mathbb Z
$$
whereby the group of deck transformations is $\mathbb Z \oplus \mathbb Z$. 
Turning to Artin's theorem, we denote by $F_k$ the free 
group on $k$ generators and $C_k$ its commutator subgroup. 
\paragraph*{Theorem 5:} (i) The covering $p$ given by (2) is regular and $p_*(\pi_1(G))$ 
is the commutator subgroup of $F_2$, the fundamental 
group of $S^1\vee S^1$. 

(ii) $p_*(\pi_1(G))$ is not finitely generated. 
\paragraph*{Proof:}  (i) Since the group of deck transformations 
of the covering is 
$\mathbb Z \oplus \mathbb Z$ and acts transitively on the fibers, the covering is regular and 
$$
\pi_1(S^1\vee S^1)/p_*(\pi_1(G)) = F_2/p_*(\pi_1(G)) = \mathbb Z\oplus \mathbb Z.\eqno(3)
$$
By virtue of (3) we conclude that $C_2 \subset p_*(\pi_1(G))$. However since $F_2/C_2$ is 
also $\mathbb Z\oplus \mathbb Z$, we get by the third isomorphism theorem, $p_*(\pi_1(G)) = C_2$.  
Indeed, since $p_*$ is injective, it establishes an isomorphism between $\pi_1(G)$ and $C_2$. 

(ii) Let $G_n = G \cap ([0, n]\times [0, n])$,  the truncated grid,  whose fundamental group is $F_{n^2}$. 
There is an obvious  retraction from $G$ onto $G_n$ and so a surjection 
$r_* : \pi_1(G) \longrightarrow \pi_1(G_n)$. If $\pi_1(G)$ were finitely generated with say $k$ generators, 
we would have a  surjection $s : F_k \longrightarrow \pi_1(G)$ thereby providing  a surjection 
$$
r_* \circ s : F_k \longrightarrow F_{n^2}. 
$$ 
On abelianizing, we get an epimorphism $A(F_k) \longrightarrow A(F_{n^2})$ which, in view of theorem 
2(iii), is contradiction since $n$ is arbitrary.  \hfill$\square$

The first homology group $H_1(G)$ is the free abelian group with countably infinite rank and $H_1(S^1 \vee S^1)$ is the 
free abelian group of rank two. 
The induced map 
$H_1(p) : H_1(G) \longrightarrow H_1(S^1 \vee S^1)$  
is trivial in stark contrast with $p_*$ which is a monomorphism.  
\paragraph*{Corollary 6:} The induced map in homology $H_1(p) : H_1(G) \longrightarrow H_1(S^1 \vee S^1)$ is trivial.
\paragraph*{Proof:} Recall that the Poincar\'e-Hurewicz map $\Pi_X  : \pi_1(X, x_0) \longrightarrow H_1(X)$ 
is a surjective group homomorphism with kernel as the commutator subgroup $[\pi_1(X, x_0), \pi_1(X, x_0)]$. 
Further, $\Pi_X$ is  
natural in the sense that the following diagram commutes: 
\begin{center}
$
\begin{CD}
\pi_1(G, x_0) @> {p}_{*} >>  \pi_1(S^1 \vee S^1, (1, 1)) \\
 @V{\Pi_G}VV @VV{\Pi_{S^1 \vee S^1}} V\\
 H_1(G) @> H_1(p) >> H_1(S^1 \vee S^1) \\
\end{CD}
$
\end{center}
Since by theorem 5 the image of $p_*$ is precisely the commutator subgroup, the composite 
$\Pi_{S^1 \vee S^1} \circ p_*$ is the zero map whereby $H_1(p) \circ \Pi_G$ is the zero map and $\Pi_G$ being surjective, 
we infer $H_1(p)$ vanishes. \hfill$\square$
\section{Connections with complex analysis:} The definition of basis of a free group mirrored the notion of basis of a vector space but the 
analogy quickly breaks down. If $W$ is a vector subspace of $V$, the dimension of $W$ does not exceed the dimension of $V$ but 
the corresponding result for free groups fails rather spectacularly as Artin's theorem shows. Massey \cite{massey} (p. 202)
has illustrated via covering spaces that the free group on seven generators embeds in the free group on two generators. 
We shall provide here an example of a covering projection which shows that the free group on four generators embeds in the free 
group on two generators. This example incidentally is an irregular covering and comes from the theory of 
Riemann surfaces \cite{reyssat} and \cite{forster} (p. 39). 
\paragraph*{Theorem 7:} The map $p : \mathbb C - \{\pm 1, \pm 2\} \longrightarrow \mathbb C - \{\pm 2\}$ 
given by 
$$
p(w) = w^3 - 3w
$$
is a three sheeted 
irregular covering projection. 
\paragraph*{Proof:} The derivative $p^{\prime}(w) = 3(w^2-1)$ vanishes precisely at $w = \pm 1$ and 
so the inverse function theorem says that $p$ is a local homeomorphism  
on $\mathbb C - \{\pm 1\}$. Now $p(1) = -2$ and $p(-1) = 2$ and so 
the equation $w^3 - 3w = z$ has a double root precisely when $z = \pm 2$ and three distinct roots for all other values of 
$z$. The removal of the points $z = \pm 2$ from the target space and the corresponding removal of the four points 
$p^{-1}\{\pm 2\}$ 
from the domain now ensures that the cardinality of the fibers 
$p^{-1}(z)$ is three throughout  $z \in \mathbb C - \{\pm 2\}$.   
The map $p : \mathbb C - \{\pm 1, \pm 2\} \longrightarrow \mathbb C - \{\pm 2\}$ 
is also surjective and a proper map and in fact it is a 
covering projection \cite{lima} (pp. 127-128). 
Let us determine its group of deck transformations. Each deck-transformation being a lift of the 
holomorphic function $p$, is itself a holomorphic function 
$\phi : \mathbb C - \{\pm 1, \pm 2\} \longrightarrow \mathbb C - \{\pm 1, \pm 2\}$ satisfying 
$$
(\phi(z))^3 - 3 \phi(z) = z^3 - 3z. \eqno(4)
$$
We see that $\phi(z)$ must be bounded in a neighborhood of the punctured points $\pm 1, \pm 2$. Riemann's removable 
singularities theorem now implies that  $\phi(z)$ extends as an entire function. One can now 
factorize\footnotemark\footnotetext{This is possible only because of the fortuitous circumstance of 
obtaining a  quadratic for $\phi(z)$} equation (4) 
to conclude that 
$\phi(z) = z$ is the only holomorphic possibility. A more interesting approach would be to observe that 
$|\phi(z)| = O(|z|)$ for large $|z|$ 
whereby $\phi(z) = a + bz$. The values of 
$a$ and $b$ are easily determined as being $0$ and $1$ respectively. 
In any case the group of deck-transformations is trivial and by virtue of theorem 4(ii), 
the covering is irregular.  \hfill$\square$
\paragraph*{Theorem 8 (Nielsen-Schreier):} If $G$ is a free group then all 
its subgroups are also free. If $G$ has rank $k$ 
and $H$ is a subgroup of $G$ of rank $l$ then,   
$$
l = (k-1)[G : H] + 1 \eqno\square
$$
We mention here that Jakob Nielsen proved the first part of the theorem \cite{chandler} (p. 84) when the subgroup is 
finitely generated\footnotemark\footnotetext{The authors of \cite{chandler} 
mention on p. 84 that prior to Nielsen, Max Dehn gave the 
topological proof of the fact that a subgroup of a free group is free.}  
and was generalized by Otto Schreier who also gave the formula stated in the theorem 
\cite{chandler} (pp. 96-97).  
The proof via covering spaces is available on p. 204 of \cite{massey}. 
If we have  a covering projection $p : ({\tilde X}, {\tilde x}_0) \longrightarrow (X, x_0)$ 
for which $\pi_1(X, x_0)$ is a free group, then we may  take $G = \pi_1(X, x_0)$,  
$H = p_*(\pi_1({\tilde X}, {\tilde x}_0))$ and the index $[G : H]$ is then the cardinality of the fiber  
$p^{-1}(x_0)$. Recalling that 
the fundamental group of the plane punctured at $k$ points is precisely $F_k$, we have  $k = 2$ and $l = 4$ 
for the covering projection of theorem 7. Nielsen-Schreier formula  gives the value of the 
index $[G : H]$ as three which is precisely 
the cardinality of the fibers.   One can of course take up more exotic examples such as $p(w) = w^5 - 5w$ in lieu of the one 
in theorem 7 for which the numbers are $k = 4$, $l = 16$ and $[G : H] = 5$.

%Let us now consider the polynomial $p(w) = w^5 - 5w$ whose derivative vanishes at the four points $\pm 1, \pm i$ and the 
%corresponding values under $p$ being $S = \{\mp 4, \mp 4i\}$. As in theorem 7, we get a covering projection 
%$p : \mathbb C - p^{-1}(S) \longrightarrow \mathbb C - S$ corresponding to  which we have 
For the example in Massey's book (p. 202) the Nielsen-Schreier  
formula predicts that the covering must be six sheeted which is indeed the case. 
For the infinite grid we have 
$[G : H] = \infty$ and the Nielsen-Schreier formula predicts that the subgroup $H$ is not finitely generated.

We conclude this note by drawing attention to the close analogy between 
the Nielsen-Schreier formula and the Riemann-Hurwitz formula for unbranched finite coverings of a 
compact Riemann surface \cite{forster} (p. 140). These connections are developed in  \cite{lyndon} (p. 133 ff.).  
\section*{Acknowledgement:} The author would like to thank C. S. Aravinda  from the TIFR Center for Applicable Mathematics 
and the referee  for their valuable feedback and 
suggestions towards a major improvement of the quality of this note.  

\end{document}